\documentclass[10pt,twoside,reqno]{amsart}
\usepackage{amssymb}
\calclayout

\numberwithin{equation}{section}
\makeatletter

\renewcommand{\@secnumfont}{\bfseries}

\renewcommand{\section}{\@startsection{section}{1}%
  {0mm}{.7\linespacing\@plus\linespacing}{.5\linespacing}
  {\normalfont\bfseries\centering}}

\newcommand{\bibsection}{\@startsection{section}{1}%
  {0mm}{.7\linespacing\@plus\linespacing}{.5\linespacing}
  {\normalfont\scshape\centering}}

\renewcommand{\@biblabel}[1]{#1.}

\newtheorem{thm}{\bf Theorem}[section]
\newtheorem{lem}[thm]{\bf Lemma}

\begin{document}

\vspace{1.3cm}

\title
        {A note on the new $q$-extension of Frobenius-Euler numbers and polynomials}

\author{ Taekyun Kim}

\thanks{\scriptsize }

\address{Department of  Mathematics, Kwangwoon  University, Seoul 139-701, Republic of Korea}
\email{tkkim@kw.ac.kr}

\keywords{$q-$umbral algebra,  $q-$umbral calculus, $n-$th ordinary Frobenius-Euler polynomials}
\subjclass{11B68, 11S80}

\maketitle

\begin{abstract} In this paper, we consider the new $q$-extension of Frobenius-Euler numbers and polynomials and we derive some interesting identities from the orthogonality type properties for the new q-extension of Frobenius-Euler polynomials. Finally we suggest one open question related to orthogonality of polynomials.
\end{abstract}

\pagestyle{myheadings} \markboth{\centerline{\scriptsize T. Kim}}
          {\centerline{\scriptsize A note on the new $q$-extension of Frobenius-Euler numbers and polynomials}}
\bigskip
\bigskip
\medskip
\section{\bf Introduction}
\bigskip
\medskip

For $q\in [0,1]$, we define the $q-$shifted factorials by $(a:q)_0=1$, $(a:q)_n=\prod_{i=0}^{n-1}(1-aq^i)$, $(a:q)_{\infty}=\prod_{i=0}^{\infty}(1-aq^i)$.
If $x$ is classical object, such as a complex number, its $q-$version is defined by $[x]_q={1-q^x\over 1-q}$.
As is well known, the $q-$exponential function is given by\
\begin{equation}
e_q(z)=\sum_{ n=0}^{\infty}  {z^n\over [n]_{q}!}={1\over {((1-q)z:q)_{\infty}}}, \ \ (\textnormal{see}\ [5],[7]),
\end{equation}
where $z\in \mathbb C$ with $\mid z \mid <1 .$
The $q-$derivative is defined by

\begin{equation}
D_qf(x)={f(x)-f(qx)\over (1-q)x}, \ \ (\textnormal{see}\ [5],[8],[9],[10]).
\end{equation}

Note that $\lim_{q\rightarrow 1} D_qf(x)={df(x)\over dx}$.\

The definite $q-$integral is given by
\begin{equation}
\int_{0}^{x} f(t)d_q t=(1-q) \sum_{a=0}^{\infty}f(q^a x)xq^a, \ \
(\textnormal{see}\ [4],[7],[10]).
\end{equation}

From (1.2) to (1.3), we have

\begin{equation}
\int_{0}^{x} t^n d_q t={1\over [n+1]_q}x^{n+1}, \ \ \ D_qx^n=[n]_qx^{n-1}.
\end{equation}

The $q-$Bernoulli polynomials of Hegazi and Mansour are defined by the generating function to be
\begin{equation}
\sum_{n=0}^{\infty}B_{n,q}(x) {t^n\over[n]_q!}={t\over e_q(t)-1}e_q(xt), \ \ (\textnormal{see}\ [4]).
\end{equation}

In the special case, $x=0, B_{n,q}(0)=B_{n,q}$ are called the $n-$th $q-$Bernoulli numbers.\\
By(1.5), we get
\begin{equation}\begin{split}
 \sum_{n=o}^{\infty}B_{n,q}(x){t^n\over [n]_q!} &=\left(\sum_{l=0}^{\infty}B_{l,q}{t^l\over [l]_q!}\right) \left( \sum_{m=0}^{\infty} {x^m\over [m]_q!}t^m\right)\\
&=\sum_{n=0}^{\infty}\left(\sum_{l=0}^{n}  \binom{n}{l}_q B_{l,q}x^{n-l}\right) {t^n\over [n]_q!}
\end{split} \end{equation}
where $\binom{n}{k}_q={[n]_q!\over [k]_q![n-k]_q!}, \ \ [n]_q!=[n]_q[n-1]_q\cdot\cdot\cdot [2]_q[1]_q.$

Thus, from (1.6), we have
\begin{equation}
B_{n,q}(x)=\sum_{l=0}^n \binom{n}{l}_q B_{l,q}x^{n-l}.
\end{equation}

In [1.6], the new $q-$extension of Euler polynomials are also defined by the generating function to be
\begin{equation}
{2\over e_q(t)+1}e_q(xt)=
\sum_{n=0}^{\infty}E_{n,q}(x){t^n\over[n]_q!}.
\end{equation}

In the special case, $x=0, \ E_{n,q}(0)=E_{n,q}$ are called the $n-$th $q-$Euler numbers.\\
By (1.8), we easily get
\begin{equation}
E_{n,q}(x)=\sum_{k=0}^n\binom{n}{k}_q E_{k,q}x^{n-k}, \ \ (\textnormal{see}\ [7]).
\end{equation}

More than five decades age, Carlitz(\textnormal{see}\ [8],[9]) defined a $q-$extension of Bernoulli and Euler polynomials.
In a recent paper(\textnormal{see}\ [10]), B.A. Kupershmidt constructed reflection symmetries of $q-$Bernoulli numbers and polynomials. From the methods of B.A. Kupershmidt, Hegazi and Mansour derived some interesting identities and properties related to $q-$Bernoulli and Euler polynomials. Recently, several authors have studied various $q-$extention of Bernoulli, Euler and Genocchi polynomials(\textnormal{see}\ [1],[2],[4]-[14]). Let $\mathcal{F}$ be the set of all formal power series in variable $t$ over $\mathbb{C}$ with\\
\begin{equation}
\mathcal{F}=\left\{f(t)=\sum_{k=0}^{\infty} {a_k\over [k]_q}t^k
\mid a_k \in\mathbb{C}\right\}.
\end{equation}

Let $\mathbb{P}=\mathbb{C}[x]$ and let $\mathbb{P^*}$ be the
vector space of all linear functionals on $\mathbb{P}$. $<L \mid
p(x)>$ denotes the action of linear functional $L$ on the
polynomials $p(x)$, and it is well know that the vector space
operation on $\mathbb{P^*}$ defined by $< L+ M \mid  p(x)>=< L
\mid p(x)>+< M\mid p(x)>, \ < cL \mid p(x)>= c< L\mid p(x)>$,
where $c$ is complex constant. For
$f(t)=\sum_{k=0}^{\infty}={a_k\over [k]_q!}t^k\in {\mathcal F},$
Let us define the linear functional on $\mathbb{P}$ by setting
\begin{equation}
<f(t)\mid x^n>=a_n\ \ (n\geq 0).
\end{equation}

By (1.11), we get
\begin{equation}
<t^k\mid x^n>=[n]_q!\ \delta_{n,k},\ \ (n, k \geq 0),
\end{equation}
where $\delta_{n,k}$ is Kronecker's symbol.

Let $f_L(t)=\sum_{k=0}^{\infty} <L\mid x^k>{{t^k}\over [k]_q!}.$
Then, by (1.12), we get
\begin{equation}
<f_L(t)\mid x^n>=<L\mid x^n>.
\end{equation}

Additionally, the map $L\rightarrow f_L(t)$ is a vector space
isomorphism from $\mathbb P^*$ onto $\mathbb F$. Henceforth,
$\mathbb F$ denotes both the
algebra of formal power series in $t$ and the vector space of all linear functionals on $\mathbb P$, and so an element $f(t)$ of $\mathbb F$ will be thought as a formal power series and a linear functional. We call $\mathbb F$ the $q-$umbral algebra. The $q-$umbral calculus is the study of $q-$umbral algebra. By (1.1) and (1.11), we easily see that $<e_q(yt)\mid x^n>=y^n$ and so $<e_q(yt)\mid p(x)>=p(y)$ for $p(x)\in \mathbb P$. The order $O(f(t))$ of the power series $f(t)\neq 0$ is the smallest positive integer for which $a_k$ does not vanish. If $O(f(t))=0$, then $f(t)$ is called an invertible series. If $O(f(t))=1$, then $f(t)$ is called a delta series(\textnormal{see}[13], [14]).\\

For $f(t)\in \mathbb F(t),\  p(x)\in \mathbb P$, we have
\begin{equation}
f(t)=\sum_{k=0}^{\infty}<f(t)\mid t^k> {t^k\over [k]_q!},\ \ p(x)=\sum_{k=0}^{\infty}<t^k\mid p(x)>{x^k\over [k]_q!}.
\end{equation}

Thus, by (1.14), we get
\begin{equation}
p^{(k)}(0)=<t^k\mid p(x)>,\ \ <1\mid p^{(k)}(x)>=p^{(k)}(0),
\end{equation}
where $p^{(k)}(x)=D_q^kp(x)\ (\textnormal{see}[13], [14]).$

From (1.15), we have
\begin{equation}
t^kp(x)=p^{(k)}(x)=D_q^kp(x).
\end{equation}

For $f(t), g(t)\in \mathbb F$ with $O(f(t))=1,\ O(g(t))=0$, there exists a unique sequence $S_n(x)\ (deg S_n(x)=n)$ of polynomial such that $<g(t)f(t)^k \mid S_n(x)>=[n]_q!\ \delta_{n,k},\ (n,\ k\geq 0)$. The sequence $S_n(x)$ is called the $q-$Sheffer sequence for $(g(t), f(t))$ which is denoted by $S_n(x)\sim (g(t), f(t)).$ \\
For $S_n(x)\sim (g(t), f(t))$,we have
\begin{equation}
{1\over g(\overline{f}(t))}e_q(y\overline{f}(t))=\sum_{k=0}^{\infty}{S_k(y)\over [k]_q!}t^k,\ for \ \ all\  y \in \mathbb{C},
\end{equation}
where $\overline{f}(t)$ is the compositional inverse of $f(t)$ (\textnormal{see}[13]).\\
 Many researchers told that the properties of the $q-$extension of Frobenius-Euler polynomials are valuable and worthwhile in the areas of both number theory and mathematical physics (\textnormal{see}[1],[2],[4]-[12]). In this paper, we consider new approach to $q-$extension of Frobenius-Euler numbers and polynomials which are derived from $q-$umbral Calculus and we give some interesting identities of our $q-$Frobenius-Euler numbers and polynomials. Finally, we suggest one open question related to orthogonality of Carlitz's $q-$Bernoulli polynomials.

\section{\bf  $q$-extension of Frobenius-Euler polynomials}
\bigskip
\medskip

In this section, we assume that $\lambda \in \mathbb C$ with $\lambda \neq 1$. Now, we consider new $q-$extension of Frobenius-Euler polynomials which are derived form the generating function as follows:
\begin{equation}
{1-\lambda \over e_q(t)-\lambda}e_q(xt)= \sum_{n=0}^{\infty}H_{n,q}(x\mid\lambda){t^n\over [n]_q!}.
\end{equation}

In the special case, $x=0,\ H_{n,q}(0\mid \lambda)= H_{n,q}(\lambda)$ are called the $n-$th $q-$Frobenius-Euler numbers. Note that $  \lim_{q\rightarrow 1}  H_{n,q}(x\mid \lambda)=H_{n}(x\mid \lambda)$, where $H_{n}(x\mid \lambda)$ are the $n-$th ordinary Frobenius-Euler polynomials.\\
By (2.1), we easily see that
\begin{equation}\begin{split}
\sum_{n=0}^{\infty}H_{n,q}(x\mid \lambda) {t^n\over [n]_q!}&= \left(\sum_{l=0}^{\infty}H_{l,q}(\lambda){t^l\over [l]_{q}!} \right) \left(\sum_{m=0}^{\infty}{x^m\over [m]_{q}!}t^m \right)\\
&=\sum_{n=0}^{\infty}\left(\sum_{l=0}^{n}\ \binom{n}{l}_q H_{l,q}(\lambda)x^{n-l}\right){t^n\over [n]_q!}.
\end{split}\end{equation}

Thus, from (2.2), we have

\begin{equation}
H_{n.q}(x\mid \lambda)=\sum_{l=0}^{n} \binom{n}{l}_q H_{l,q}(\lambda)x^{n-l}.
\end{equation}

By the definition of $q-$Frobenius-Euler numbers, we get

\begin{equation}\begin{split}
1-\lambda &=\left(\sum_{l=0}^{\infty}H_{l,q}(\lambda) {t^l\over [l]_q!} \right)\left(e_q(t)-\lambda \right)\\
&=\sum_{n=0}^{\infty}\left\{\sum_{m=o}^n\binom{n}{m}_q H_{n-m,q}(\lambda) \right\}{t^n\over[n]_q!}-\lambda \sum_{n=0}^{\infty}H_{n,q}(\lambda){t^n\over[n]_q!}\\
&=\sum_{n=0}^{\infty}\left\{ H_{n,q}(1 \mid\lambda)- \lambda H_{n,q}(\lambda)\right\}{t^n\over [n]_q!}.
\end{split}\end{equation}

By comparing the coefficients on the both sides, we get
\begin{equation}
H_{0,q}(\lambda)=1,\ H_{n,q}(1\mid \lambda)- \lambda H_{n,q}(\lambda)=
\left\{
\begin{array}{cc}
\ 1- \lambda , & \text{if }n=0 \\
0, & \text{if\ }n> 0\text{.}%
\end{array}%
\right.
\end{equation}

Therefore, by (2.3) and (2.5), we obtain the following theorem.

\begin{thm}\label{THEOREM 2.1} For $n \geq 0$, we have
\begin{equation*}
H_{n,q}(x\mid \lambda)=\sum_{l=0}^{n}\binom{n}{l}_q H_{n-l,q}(\lambda)x^l .\\
\end{equation*}
Moreover,\\
\begin{equation*}
H_{0,q}(\lambda)=1,\ H_{n,q}(1\mid \lambda)-\lambda H_{n,q}(\lambda)= (1-\lambda)\delta_{0,n}.
\end{equation*}
\end{thm}

\medskip

For examples, $H_{0,q}(\lambda)=1,\ H_{1,q}(\lambda)={-1\over
1-\lambda}, H_{2,q}(\lambda)={\lambda +q\over
(1-\lambda)^2},\cdots. $

From (1.17) and (2.1), we have
\begin{equation}
 H_{n,q}(x\mid \lambda) \sim \left({e_q(t)-\lambda\over 1-\lambda},\ t
 \right),
\end{equation}

and
\begin{equation}
{1-\lambda \over e_q(t)-\lambda}e_q(xt) =\sum_{n=0}^{\infty}H_{n,q}(x\mid \lambda){t^n\over [n]_q!}. \\
\end{equation}

Thus, by (2.7), we get

\begin{equation}
{1-\lambda \over e_q(t)-\lambda}x^n =H_{n,q}(x\mid \lambda),\ (n\geq 0). \\
\end{equation}

From (1.4), (1.16) and (2.8), we have
\begin{equation}\begin{split}
tH_{n,q}(x\mid \lambda)&={1-\lambda \over e_q(t)-\lambda}tx^n =[n]_q{1-\lambda \over e_q(t)-\lambda}x^{n-1}\\
&=[n]_q H_{n-1,q}(x\mid \lambda).
\end{split}\end{equation}

Therefore, by (2.8) and (2.9), we obtain the following lemma.

\begin{lem}\label{Lemma 2.2} For $n\geq 0$, we have
\begin{equation*}
H_{n,q}(x\mid \lambda)={1-\lambda \over e_q(t)-\lambda}x^n, \
tH_{n,q}(x\mid \lambda)=[n]_q H_{n-1,q}(x\mid \lambda).
\end{equation*}
\end{lem}

By (2.9), we get
\begin{equation}\begin{split}
 \left <{e_q(t)-\lambda    \over  1-\lambda }t^k \mid  H_{n,q}(x\mid \lambda)\right> &={[k]_{q}!\over 1-\lambda} \binom {n}{k}_q \left<e_q(t)-\lambda \mid H_{n-k,q}(x\mid\lambda) \right>\\
  &={[k]_{q}!\over 1-\lambda} \binom {n}{k}_q  \left\{H_{n-k,q}(1\mid\lambda) - \lambda H_{n-k,q}(\lambda)
  \right\}.
\end{split}\end{equation}

From (2.6), we have
\begin{equation}
 \left <{e_q(t)-\lambda    \over  1-\lambda }t^k \mid  H_{n,q}(x\mid \lambda)\right>=[n]_{q}! \
 \delta_{n,k}.
\end{equation}

  By (2.10) and (2.11), we get

\begin{equation}\begin{split}
0&={1\over 1-\lambda}\left(  H_{n-k,q}(1\mid \lambda)-\lambda H_{n-k,q}(\lambda) \right)\\
&=\sum_{l=0}^{n-k-l}\binom {n-k}{l}_q
H_{l,q}(\lambda)+(1-\lambda)H_{n-k,q}(\lambda),
\end{split}\end{equation}
  where $n, k \in \mathbb{Z}_{\geq 0} $ with $n>k$.

 Thus, from (2.12), we have

\begin{equation}
 H_{n-k,q}(\lambda)={1\over \lambda -1}\sum_{l=0}^{n-k-l}\binom {n-k}{l}
 H_{l,q}(\lambda),
\end{equation}
  where $n, k \in \mathbb{Z}_{\geq 0}$ with $n>k$.

  Therefore, by (2.13), we obtain the following theorem.

 \begin{thm}\label{THEOREM 2.3} For $n \geq 1$, we have
\begin{equation*}
 H_{n,q}(\lambda)={1\over \lambda -1}\sum_{l=0}^{n-l}\binom {n}{l}
 H_{l,q}(\lambda).
\end{equation*}
\end{thm}
\medskip

 From (1.4) and (2.3), we have
 \begin{equation}\begin{split}
 \int_x^{x+y}H_{n,q}(u \mid \lambda)d_qu &={1\over [n+1]_q}\sum_{l=0}^{n}\binom {n+1}{l+1}_q H_{n-l,q}(\lambda) \left\{ (x+y)^{l+1}-x^{l+1} \right\}\\
 &={1\over [n+1]_q}\sum_{l=0}^{n+1}\binom {n+1}{l}_q H_{n+1-l,q}(\lambda) \left\{ (x+y)^{l}-x^{l} \right\}\\
 &={1\over [n+1]_q} \left\{H_{n+1,q}(x+y\mid \lambda) -H_{n+1,q}(x\mid \lambda)\right\}.\\
\end{split}\end{equation}

Thus, by (2.14), we get
\begin{equation}\begin{split}
 \left<{e_q(t)-1\over t} \mid H_{n,q}(x\mid \lambda)\right> &= {1\over [n+1]_q} \left<{e_q(t)-1\over t}\mid t H_{n+1,q}(x\mid \lambda)\right>\\
 &= {1\over [n+1]_q} \left<e_q(t)-1\mid H_{n+1,q}(x\mid \lambda)\right>\\
  &= {1\over [n+1]_q} \left\{ H_{n+1,q}(1\mid \lambda) - H_{n+1,q}( \lambda)\right\}\\
  &=\int_0^{1}H_{n,q}(u \mid \lambda)d_qu.
\end{split}\end{equation}

 Therefore, by (2.15), we obtain the following theorem.

 \begin{thm}\label{THEOREM 2.4} For $n \geq 0$, we have
\begin{equation*}
\left<{e_q(t)-1\over t} \mid H_{n,q}(x\mid
\lambda)\right>=\int_0^{1}H_{n,q}(u \mid \lambda)d_qu.
\end{equation*}
\end{thm}
\medskip

Let
\begin{equation}
\mathbb{P}_n=\left\{p(x)\in\mathbb{C}[x]\mid deg\ p(x)\leq n
\right\},
\end{equation}
be the $(n+1)-$dimensional vector space over $\mathbb{C}$. For
$p(x)\in \mathbb{P}_n$, let us assume that

\begin{equation}
 p(x)=\sum_{k=0}^{n} C_k H_{k,q}(x\mid \lambda), \ \ (n\geq 0).
\end{equation}

Then, by (2.6), (2.11) and (2.17), we get

\begin{equation}\begin{split}
 \left<\left({e_q(t)-\lambda\over 1-\lambda}\right)t^k \mid p(x)\right> &= \sum_{l=0}^{n} C_l  \left<\left({e_q(t)-\lambda\over 1-\lambda}\right)t^k\mid  H_{l,q}(x\mid \lambda)\right>\\
 &= \sum_{l=0}^{n} C_l[l]_q! \ \delta_{l,k}=[k]_q! \ C_k.
\end{split}\end{equation}

By (2.18), we get
\begin{equation}\begin{split}
 C_k&={1\over [k]_q!}\left<\left({e_q(t)-\lambda\over 1-\lambda}\right)t^k \mid p(x)\right>\\
 &= {1\over (1-\lambda)[k]_q!} \left\{p^{(k)}(1)-\lambda p^{(k)}(0)  \right\},
\end{split}\end{equation}
where $p^{(k)}(x)=D_q^kp(x)$.
Therefore, by (2.17) and (2.19), we obtain the following theorem.

\begin{thm}\label{THEOREM 2.5} For $p(x)\in \mathbb{P}_n$, let
$p(x)=\sum_{k=0}^{n} C_k H_{k,q}(x\mid \lambda)$\\
Then we have
\begin{equation*}\begin{split}
C_k={1\over [k]_q! \ (1-\lambda)} \left\{p^{(k)}(1)-\lambda p^{(k)}(0)  \right\},
\end{split}\end{equation*}
where $p^{(k)}(x)=D_q^kp(x)$.
\end{thm}
\medskip

Let us take $p(x)=B_{n,q}(x)$ with

\begin{equation}
B_{n,q}(x)=p(x)=\sum_{k=0}^{n} C_k H_{k,q}(x\mid \lambda).
\end{equation}

Then, by Theorem 2.5, we get
\begin{equation}\begin{split}
C_k&={1\over [k]_q! \ (1-\lambda)}\\
&\times \left\{ [n]_q [n-1]_q \cdot\cdot\cdot [n-k+1]_q B_{n-k,q}(1)-\lambda [n]_q\cdot\cdot\cdot [n-k+1]_qB_{n-k,q}\right\}\\
&= \binom {n}{k}_q {1\over 1-\lambda} \left\{
B_{n-k,q}(1)-\lambda B_{n-k,q}\right\}.
\end{split}\end{equation}

Therefore, by (2.21), we obtain the following theorem.

\begin{thm}\label{Theorem 2.6}   For $n, \ k\geq 0$, with $n-k>0 $, we have
\begin{equation*}
 B_{n,q}(x)={1\over 1-\lambda}\sum_{k=0}^{n} \binom {n}{k}_q \left\{  B_{n-k,q}(1)-\lambda B_{n-k,q}\right\}
 H_{k,q}(x\mid \lambda).
\end{equation*}
\end{thm}
\medskip

Let us take $p(x)=x^n$ with
\begin{equation}
x^n=p(x)=\sum_{k=0}^n C_kH_{k,q}(x\mid \lambda).
\end{equation}

From (1.5), we have

\begin{equation}
C_k={1\over 1-\lambda} \binom{n}{k}_q -{\lambda\over 1-\lambda}\binom{n}{k}_q 0^{n-k}.\\
\end{equation}

Thus, we have
\begin{equation*}\begin{split}
x^n &={1\over 1-\lambda} \sum_{k=0}^{n}\binom{n}{k}_q H_{k,q}(x\mid \lambda)-{\lambda\over 1-\lambda}H_{n,q}(x\mid \lambda)\\
&={1\over 1-\lambda} \sum_{k=0}^{n-1}\binom{n}{k}_q H_{k,q}(x\mid \lambda)-H_{n,q}(x\mid \lambda).\\
\end{split}\end{equation*}

By the same method, we easily see that

\begin{equation*}
E_{n,q}(x)={1\over 1-\lambda} \sum_{k=0}^{n}\binom{n}{k}  \left\{ E_{n-k,q}(1)-\lambda E_{n-k,q}\right\}
  H_{k,q}(x\mid \lambda).
\end{equation*}

The $q-$extension of Frobenius-Euler polynomials of order $r \
(r\in \mathbb{N})$ are defined by the generating function to be

\begin{equation}
\left({1-\lambda \over e_q(t)-\lambda} \right)^r e_q(xt)=
\sum_{n=0}^{\infty}H_{n,q}^{(r)}(x\mid \lambda) {t^n\over [n]_q!}.
\end{equation}

In the special case, $x=0, \ H_{n,q}^{(r)}(0\mid\lambda)=H_{n,q}^{(r)}(\lambda)$ are called the $n-$th $q-$Frobenius-Euler numbers of order $r$.

From (2.24), we can derive the following equation (2.25):
\begin{equation}
H_{n,q}^{(r)}(x\mid \lambda)= \sum_{l=0}^{n} \binom{n}{l}_q   H_{n-l,q}^{(r)}(\lambda)x^l.
\end{equation}
by (1.17) and (2.24), we get

\begin{equation}
H_{n,q}^{(r)}(x\mid \lambda) \sim \left( \left(
{e_q(t)-\lambda\over 1-\lambda}\right)^r, \ t \right).
\end{equation}

Thus, from (2.26), we have

\begin{equation}
      \left<\left( {e_q(t)-\lambda\over 1-\lambda}\right)^r t^k \mid H_{n,q}^{(r)}(x\mid \lambda )
       \right>  =[n] _q! \ \delta_{n,k}.
\end{equation}

For $p(x) \in\mathbb{P}_n$, let us assume that
\begin{equation}
p(x)=\sum_{k=0}^{n}C_k^{r} H_{k,q}^{(r)} (x\mid \lambda),\  (n\geq
0).
\end{equation}

From (2.27) and (2.28), we can derive the following equation (2.29):

\begin{equation}\begin{split}
     \left<\left( {e_q(t)-\lambda\over 1-\lambda}\right)^r t^k \mid  p(x)  \right>
     &=\sum_{l=0}^{n}C_l^{r} \left<\left( {e_q(t)-\lambda\over 1-\lambda}\right)^r t^k \mid H_{l,q}^{(r)}(x\mid \lambda  \right>\\
      &=\sum_{l=0}^{n} C_l^{r} [l]_q!  \delta_{l,k}= C_k^{r} [k]_q!.
\end{split}\end{equation}

Thus, by (2.29), we get

\begin{equation}\begin{split}
  C_k^r&={1\over[k]_q!}\left<\left( {e_q(t)-\lambda\over 1-\lambda}\right)^r t^k \mid  p(x)  \right> \\
   &={1\over [k]_q!(1-\lambda)^r}\sum_{j=0}^r\binom{r}{j}(-\lambda)^{r-j} \left< \left( e_q(t) \right)^j \mid t^kp(x)  \right>\\
   &={1\over [k]_q!(1-\lambda)^r}\sum_{j=0}^r \binom{r}{j}(-\lambda)^{r-j}\sum_{l=0}^{\infty}\sum_{l_1+\cdot\cdot\cdot+l_j=l}{\left< 1\mid p^{(k+l)}(x) \right> \over [l_1]_q!\cdot\cdot\cdot[l_j]_q!}\\
   &={1\over [k]_q!(1-\lambda)^r}\sum_{l=0}^{\infty}    \sum_{j=0}^{r}\sum_{l_1+\cdot\cdot\cdot+l_j=l}
   { \binom{r}{j}(-\lambda)^{r-j} \over
   [l_1]_q!\cdot\cdot\cdot[l_j]_q!}p^{(k+l)}(0) .
\end{split}\end{equation}

Therefore, by (2.28) and (2.30), we obtain the following theorem.

\begin{thm}\label{THEOREM 2.7}   For $p(x)\in \mathbb{P}_n $, let $p(x)=\sum_{k=0}^{n}
C_k^{r} H_{k,q}^{(r)} (x\mid \lambda).$\\\\
Then we have
\begin{equation*}\begin{split}
  C_k^r={1\over [k]_q!(1-\lambda)^r}\sum_{l=0}^{\infty}
  \sum_{j=0}^{r}\sum_{l_1+\cdot\cdot\cdot+l_j=l} { \binom{r}{j}(-\lambda)^{r-j}
  \over [l_1]_q!\cdot\cdot\cdot[l_j]_q!}p^{(k+l)}(0)  ,\\
\end{split}\end{equation*}
where $p^{(k+l)}(x)=D_q^{k+l}p(x)$.
\end{thm}
\medskip

From (1.17) and  (2.30), we have

\begin{equation}
 \sum_{n=0}^{\infty} H_{n,q}^{(r)}(x\mid \lambda){t^n\over [n]_q!}
 =\left( 1-\lambda\over e_q(t)-\lambda \right)^re_q(xt).
\end{equation}

Thus, from (2.31), we have.

\begin{equation}
  H_{n,q}^{(r)}(x\mid \lambda)=\left( 1-\lambda\over e_q(t)-\lambda \right)^rx^n, \ \ (n\geq 0),
\end{equation}

and

\begin{equation}\begin{split}
  tH_{n,q}^{(r)}(x\mid \lambda)&=\left( 1-\lambda\over e_q(t)-\lambda \right)^rtx^n=[n]_q\left( 1-\lambda\over e_q(t)-\lambda \right)^rx^{n-1}\\
  &=[n]_qH_{n-1,q}^{(r)}(x\mid \lambda).
\end{split}\end{equation}

Indeed,
\begin{equation}
 \left<\left( 1-\lambda\over e_q(t)-\lambda \right)^r e_q(yt)\mid x^n\right>
 = H_{n,q}^{(r)}(y\mid \lambda)=\sum_{l=0}^{n} \binom{n}{l}_q
 H_{n-l,q}^{(r)}(\lambda)y^l,
 \end{equation}

and
\begin{equation}\begin{split}
 \left<\left( 1-\lambda\over e_q(t)-\lambda \right)^r \mid x^n\right>&= \sum_{m=0}^{\infty}\left( \sum_{i_1+\cdot\cdot\cdot+i_r=m} {H_{i_1,q}(\lambda)\cdot\cdot\cdot H_{i_r,q}(\lambda)\over [i_1]_q!\cdot\cdot\cdot [i_r]_q!}  \right)\left< t^m\mid x^n\right>\\
 &= \sum_{i_1+\cdot\cdot\cdot+i_r=n} {[n]_q!\over [i_1]_q!\cdot\cdot\cdot [i_r]_q!}  H_{i_1,q}(\lambda)\cdot\cdot\cdot H_{i_r,q}(\lambda)\\
 &= \sum_{i_1+\cdot\cdot\cdot+i_r=n}
 \binom{n}{i_1,\cdot\cdot\cdot,i_r}_q  H_{i_1,q}(\lambda)\cdot\cdot\cdot H_{i_r,q}(\lambda),\\
 \end{split} \end{equation}
where $\binom{n}{i_1,\cdot\cdot\cdot,i_r}_q= {[n]_q!\over [i_1]_q!\cdot\cdot\cdot [i_r]_q!}.$

From (2.34), we note that
\begin{equation}
 \left<\left( 1-\lambda\over e_q(t)-\lambda \right)^r \mid x^n\right>= H_{n,q}^{(r)}(\lambda), \ (n\geq
 0).
 \end{equation}

Therefore, by (2.35) and (2.36), we obtain the following theorem.

\begin{thm}\label{THEOREM 2.8}   For $n\geq 0 $, we have
\begin{equation*}\begin{split}
 H_{n,q}^{(r)}(\lambda)=\sum_{i_1+\cdot\cdot\cdot+i_r=n} \binom{n}{i_1,\cdot\cdot\cdot,i_r}_q
 H_{i_1,q}(\lambda)\cdot\cdot\cdot H_{i_r,q}(\lambda).\\
\end{split}\end{equation*}
\end{thm}
\medskip

Let us take $p(x)=H_{n,q}^{(r)}(x\mid \lambda)\in \mathbb{P}_n$ with\\

$H_{n,q}^{(r)}(x\mid \lambda)=p(x)=\sum_{k=0}^{n}C_k H_{k,q}(x\mid \lambda)$.\\
Then, by Theorem2.5, we get
\begin{equation}\begin{split}
C_k&={1\over \left(1-\lambda \right) [k]_q!} \left<\left( e_q(t)-\lambda \right) t^k\mid p(x)\right>\\
&= {1\over 1-\lambda}\binom{n}{k}_q  \left<e_q(t)-\lambda \mid H_{n-k,q}^{(r)}(x\mid \lambda)\right>\\
&= {\binom{n}{k}_q\over 1-\lambda} \left\{ H_{n-k,q}^{(r)}(1\mid \lambda)- \lambda H_{n-k,q}^{(r)}(\lambda)\right\}.\\
\end{split} \end{equation}

From (2.24), we can readily derive the following equation (2.38):

\begin{equation}\begin{split}
&\sum_{n=0}^{\infty} \left\{ H_{n,q}^{(r)}(1\mid \lambda)- \lambda H_{n,q}(\lambda)\right\}{t^n\over [n]_q!}\\
&=\left( 1-\lambda\over e_q(t)-\lambda \right)^r  \left( e_q(t)-\lambda \right)=(1-\lambda)\left( 1-\lambda\over e_q(t)-\lambda \right)^{r-1} \\
&=(1-\lambda)\sum_{n=0}^{\infty}H_{n,q}^{(r-1)}(\lambda){t^n\over [n]_q!}.
\end{split} \end{equation}

By comparing the coefficients on the both sides of (2.38), we get

\begin{equation}
H_{n,q}^{(r)}(1\mid \lambda)-\lambda H_{n,q}^{(r)}(\lambda)=(1-\lambda)H_{n,q}^{(r-1)}(\lambda).
 \end{equation}

Therefore, by (2.37) and (2.39), we obtain the following theorem.

\begin{thm}\label{THEOREM 2.9}   For $n\in \mathbb{Z}_{\geq 0}$ and $r\in \mathbb{N} $, we have
\begin{equation*}
 H_{n,q}^{(r)}(x\mid\lambda)=\sum_{k=0}^{n} \binom{n}{k}_q  H_{n-k,q}^{(r-1)}(\lambda) H_{k,q}(x\mid\lambda).\\
\end{equation*}
\end{thm}
\medskip

Let us assume that $p(x)=H_{n,q}(x\mid\lambda)\in \mathbb{P}_n$ with

\begin{equation}
H_{n,q}(x\mid \lambda)=p(x)=\sum_{k=0}^{n}C_k^r
H_{k,q}^{(r)}(x\mid \lambda).
\end{equation}

Then, by Theorem 2.7, we get

\begin{equation}\begin{split}
C_k^r&={1\over \left( 1-\lambda\right)^r[k]_q!}\sum_{m=0}^{n-k}\sum_{l=0}^{r}\binom{r}{l}\left( -\lambda\right)^{r-l} \\
     & \ \ \ \ \ \ \ \ \ \ \ \ \ \ \times\sum_{i_1+\cdot\cdot\cdot+i_l=m} \binom{m}{i_1,\cdot\cdot\cdot,i_l}_q {[n]_q\cdot\cdot\cdot[n-m-k+1]_q\over [m]_q!}H_{n-m-k,q}(\lambda)\\
&={1\over \left( 1-\lambda\right)}\sum_{m=0}^{n-k}\sum_{l=0}^{r}\binom{r}{l}\left( -\lambda\right)^{r-l}\sum_{i_1+\cdot\cdot\cdot+i_l=m}\binom{m}{i_1,\cdot\cdot\cdot,i_l}_q\\
     & \ \ \ \ \ \ \ \ \ \ \ \ \ \times {[m+k]_q!\over [m]_q[k]_q!} \times{[n]_q\cdot\cdot\cdot[n-m-k+1]_q\over [m+k]_q!}H_{n-m-k,q}(\lambda)\\
&={1\over \left( 1-\lambda\right)^r}\sum_{m=0}^{n-k}\sum_{l=0}^{r}\sum_{i_1+\cdot\cdot\cdot+i_l=m}\binom{r}{l}\\
    &\ \ \ \ \ \ \ \ \ \ \ \ \ \times\binom{m}{i_1,\cdot\cdot\cdot,i_l}_q \binom{m+k}{m}_q \binom{n}{m+k}_q  (-\lambda)^{r-l} H_{n-m-k,q}(\lambda).
\end{split}\end{equation}

Therefore, by (2.40) and (2.41), we obtain the following theorem.

\begin{thm}\label{THEOREM 2.10}   For $n, r\geq 0 $, we have
\begin{equation*}\begin{split}
 &H_{n,q}(x\mid \lambda)={1\over \left( 1-\lambda\right)^r} \sum_{k=0}^{n}\\
 &\left\{ \sum_{m=0}^{n-k}\sum_{l=0}^{r}\sum_{i_1+\cdot\cdot\cdot+i_l=m}\binom{r}{l}(-\lambda)^{r-l}
\binom{m}{i_1,\cdot\cdot\cdot,i_l}_q \binom{m+k}{m}_q \binom{n}{m+k}_q   H_{n-m-k,q}(\lambda)    \right\}\\
&\times H_{k,q}^{(r)}(x\mid\lambda).
\end{split}\end{equation*}
\end{thm}
\medskip

Remark. By the same method, we can see that

\begin{equation*}\begin{split}
 &B_{n,q}^{(r)}(x)\\
 &={1\over \left( 1-\lambda\right)^r} \sum_{k=0}^{n}
 \left\{ \sum_{m=0}^{n-k}\sum_{l=0}^{r}\sum_{i_1+\cdot\cdot\cdot+i_l=m}\binom{r}{l}(-\lambda)^{r-l}
\binom{m}{i_1,\cdot\cdot\cdot,i_l}_q \binom{m+k}{m}_q \binom{n}{m+k}_q   B_{n-m-k,q}^{(r)} \right\}\\
&\times H_{k,q}^{(r)}(x\mid\lambda).
\end{split}\end{equation*}

As is well known, the Carltz's $q-$Bernoullli polynomials are defined by the generating function to be
\begin{equation}\begin{split}
f(t,x)&=\sum_{n=0}^{\infty}\beta_{n,q}{t^n\over n!}\\
&=-t\sum_{m=0}^{\infty}q^{2m}e^{[m+x]_qt}+(1-q)\sum_{m=0}^{\infty}q^{m}e^{[m+x]_qt}.
\end{split}\end{equation}

Thus, we note that $f(t,x)$ is an invertible series. From (2.42), we note that

\begin{equation}
\beta_{n,q}(x)={1\over
(1-q)^n}\sum_{l=0}^{n}\binom{n}{l}(-1)^lq^{lx}{l+1\over [l+1]_q},
 \text{ (see [3])}.
\end{equation}
In the special case, $x=0,\beta_{n,q}=\beta_{n,q}(0)$ are called the $n-$th Carlitz's $q-$Bernoulli numbers.\\
By (2.42) and (2.43), we readily see that

\begin{equation*}
\beta_{n,q}(x)=\sum_{l=0}^{n}\binom{n}{l} q^{lx} \beta_{l,q}
[x]_q^{n-l}\sim(f(t,x),t).
\end{equation*}

{\bf Open Problem}. Let
\begin{equation*}
\mathbb{P}_{n,q}=\left\{ p([x]_q)\in \mathbb{C} [[x]_q]| \deg
p([x]_q)\leq n \right\}
\end{equation*}

For $p([x]_q)\in \Bbb P_{n, q}$, let us assume that
$p([x]_q)=\sum_{k=0}^{n}C_{k,q} \beta_{k,q}(x)$.\\
Determine the coefficient $C_{k,q}$

\end{document}